\documentclass[12pt]{article}
\usepackage{amsthm,amsfonts,amssymb,amscd}

\begin{document}


\begin{center}
\large \bf Birationally rigid complete intersections \\ with a
singular point of high multiplicity
\end{center}\vspace{0.5cm}

\centerline{A.V.Pukhlikov}\vspace{0.5cm}

\parshape=1
3cm 10cm \noindent {\small \quad\quad\quad \quad\quad\quad\quad
\quad\quad\quad {\bf }\newline We prove the birational rigidity of
Fano complete intersecti\-ons of index 1 with a singular point of
high multiplicity, which can be close to the degree of the
variety. In particular, the groups of birational and biregular
automorphisms of these varieties are equal, and they are
non-rational. The proof is based on the techniques of the method
of maximal singularities, the generalized $4n^2$-inequality for
complete in\-tersection singularities and the technique of
hypertangent divisors.

Bibliography: 19 titles.} \vspace{1cm}

\vspace{1cm}

\noindent 14E05, 14E07\vspace{1cm}

\begin{flushright}
{\it To Yu.I.Manin on the occasion of his 80th birthday}
\end{flushright}

\section*{Introduction}

{\bf 0.1. Statement of the main result.} The aim of the present paper
is to prove the birational superrigidity of generic Fano complete
intersections of index 1 with a singular point of high multiplicity.
This is a generalization of the result of \cite{Pukh02a}, where a similar
fact was shown for hypersurfaces.\vspace{0.1cm}

Let $k\geqslant 2$ and $M\geqslant 2k+1$ be fixed integers,
${\mathbb P}={\mathbb P}^{M+k}$ the complex projective space.
Fix an ordered integral vector
$$
\underline{d}=(d_1,\dots,d_k)\in{\mathbb Z}^k_+,
$$
where $2\leqslant d_1\leqslant\dots\leqslant d_k$, satisfying the equality
$$
|\underline{d}|=d_1+\dots+d_k=M+k,
$$
and an integral vector
$$
\underline{\xi}=(\xi_1,\dots,\xi_k)\in{\mathbb Z}^k_+,
$$
where $1\leqslant\xi_i\leqslant d_i$ for all $i=1,\dots,k$. Set
$$
c_*=\sharp\{i\,|\,\xi_i=d_i,\,\,i=1,\dots,k\}.
$$
Assume that the inequalities
\begin{equation}\label{07.03.2017.1}
\sum^k_{i=1}[(d_i+1)(d_i+2)-\xi_i(\xi_i+1)]\geqslant4M+2d_k+2c_*-2k
\end{equation}
and
\begin{equation}\label{07.03.2017.2}
M\geqslant 3+\sum_{\xi_i\geqslant 2}(\xi_i+1)
\end{equation}
hold. Let
$$
{\cal H}=\prod^k_{i=1}H^0({\mathbb P},{\cal O}_{\mathbb P}(d_i))
$$
be the space of tuples of homogeneous polynomials $\underline{f}
=(f_1,\dots,f_k)$ of degrees $d_1,\dots,d_k$, respectively. Fix
a point $o\in{\mathbb P}$. The symbol
${\cal H}(\underline{\xi})$ stands for the subset of tuples
$\underline{f}\in{\cal H}$, such that:
\begin{itemize}

\item the set of common zeros
$$
V(\underline{f})=\{f_1=\dots=f_k=0\}
$$
is an irreducible reduced complete intersection of codimension $k$
in ${\mathbb P}$,

\item the subvariety $V(\underline{f})$ is non-singular outside the point
$o$,

\item $\mathop{\rm mult}_of_i=\xi_i$ for $i=1,\dots,k$.
\end{itemize}
If $\xi_i\geqslant 2$ for at least one $i\in\{1,\dots,k\}$, then
the subvariety $V(\underline{f})$ is singular at the point $o$. By
Grothendieck's theorem \cite{CL}, the variety $V(\underline{f})$
is factorial, so that $\mathop{\rm
Pic}V(\underline{f})\cong{\mathbb Z}H$, where $H$ is the class of
a hyperplane section. Note that for a general tuple
$\underline{f}\in {\cal H}(\xi)$ the unique singular point $o\in
V(\underline{f})$ is resolved by one blow up with a non-singular
exceptional divisor, the discrepancy of which is positive by the
assumption (\ref{07.03.2017.2}). Therefore, the singularity $o\in
V(\underline{f})$ is terminal and $V(\underline{f})$ is a Fano
variety of index 1 and dimension $M$.\vspace{0.1cm}

The main result of the present paper is the following claim.\vspace{0.1cm}

{\bf Theorem 0.1.} {\it There exists a non-empty Zariski open subset
$$
{\cal U}\subset{\cal H}(\underline{\xi}),
$$
such that for any tuple $\underline{f}\in{\cal U}$
the variety $V(\underline{f})$ is birationally superrigid.}\vspace{0.1cm}

We remind the definition and the main properties of birationally
superrigid varieties below in Subsection 0.3.\vspace{0.1cm}

{\bf Remark 0.1.} For a general tuple $\underline{f}\in{\cal
H}(\underline{\xi})$ the point $o\in V(\underline{f})$ has the
multiplicity $\mu=\xi_1\dots\xi_k$. The anticanonical degree of
the variety $V(\underline{f})$ (that is to say, its degree in the
projective space ${\mathbb P}$) is $d=d_1\dots d_k$. It is easy to
check that the infimum of the ratios $\mu/d$ over all tuples
$\underline{d}$, $\underline{\xi}$, satisfying the conditions
(\ref{07.03.2017.1}) and (\ref{07.03.2017.2}), is equal to 1, that
is, the multiplicity of the point $o\in V(\underline{f})$ can be
really high.\vspace{0.3cm}


{\bf 0.2. Regular complete intersections.} The Zariski open
subset ${\cal H}_{\rm reg}(\underline{\xi})$ is defined by a set
of local conditions, which must hold at every point
$p\in V(\underline{f})$. First we consider the conditions
for the singular point $o\in V=V(\underline{f})$. Set
$l=\sharp\{i\,|\,\xi_i\geqslant 2\}$ and let
$$
\underline{\mu}=(\mu_1,\dots,\mu_l)
$$
be the ordered (non-decreasing) tuple of multiplicities
$\xi_i\geqslant 2$, $i\in\{1,\dots,k\}$, that is,
$\mu_1\leqslant\mu_2\leqslant\dots\leqslant\mu_l$ and
$\mu_{\alpha}=\xi_{i_{\alpha}}$, where
$\{i_1,\dots,i_l\}=\{i\,|\,\xi_i\geqslant 2\}$. The integral vector
$\underline{\mu}$ will be called the {\it type} of the singularity
$o\in V$. Set
$$
\mu=\mu_1\dots\mu_l=\xi_1\dots\xi_k\quad\mbox{and}\quad
|\underline{\mu}|=\mu_1+\dots+\mu_l.
$$
The following condition is a natural condition of general position
for the singular point $o\in V$.\vspace{0.1cm}

(R0.1) Let $\varphi_{\mathbb P}\colon{\mathbb P}^+\to{\mathbb
P}$ be the blow up of the point $o$,
$Q_{\mathbb P}=\varphi^{-1}_{\mathbb P}(o)\cong{\mathbb P}^{M+k-1}$
the exceptional divisor, $V^+\subset{\mathbb P}^+$ the strict
transform of the complete intersection $V$, so that
$\varphi\colon V ^+\to V$ is the blow up of the point $o$ on $V$,
$Q=\varphi^{-1}(o)=V^+\cap Q_{\mathbb P}$ is the exceptional
divisor. The subvariety
$$
Q\subset Q_{\mathbb P}\cong{\mathbb P}^{M+k-1}
$$
is a non-singular complete intersection of codimension $l$ in its linear
span
$$
\langle Q\rangle\cong{\mathbb P}^{M+l-1}.
$$
Let us discuss the condition (R0.1) in more detail. Let
$(z_*)=(z_1,\dots,z_{M+k})$ be a system of affine coordinates with
the origin at the point $o$ and
$$
f_i=q_{i,\xi_i}+\dots+q_{i,d_i}
$$
the decomposition of the polynomial $f_i$ (in the non-homogeneous form) into
components, homogeneous in $z_*$. Then $(z_1:\dots:z_{M+k})$ form a set of
homogeneous coordinates on $Q_{\mathbb P}$, the linear span
$\langle Q\rangle$ is given by the system of equations
$$
q_{i,1}=0
$$
for $i\in\{1,\dots,k\}$ such that $\xi_i=1$, and the complete intersection
$Q\subset\langle Q\rangle$ is given by the system of equations
$$
q_{i,\xi_i}|_{\langle Q\rangle}=0
$$
for $i\in\{1,\dots,k\}$ such that $\xi_i\geqslant 2$. We will need one
more condition of general position for the singular point
$o\in V$. We use the coordinate notations introduced above.\vspace{0.1cm}

(R0.2) The set of homogeneous polynomials $q_{i,j}$,
$i=1,\dots,k$, $j=\xi_i,\dots,d_i$ forms a regular sequence in
${\cal O}_{o,{\mathbb P}}$, that is, the system of homogeneous
equations
$$
\{q_{i,j}=0\,|\,i=1,\dots,k,j=\xi_i,\dots,d_i\}
$$
defines a closed set of codimension
$\sum\limits^k_{i=1}(d_i-\xi_i+1)$ in ${\mathbb
C}^{M+k}_{z_*}$.\vspace{0.1cm}

Now let us consider the conditions for a non-singular point $p\in
V$, $p\neq o$. These conditions are identical to the regularity
conditions introduced in \cite{Pukh01}, see also \cite[Chapter
3]{Pukh13a}. Let $u_1,\dots,u_{M+k}$ be a system of affine
coordinates with the origin at the point $p$ and
$$
f_i=\Phi_{i,1}+\Phi_{i,2}+\dots+\Phi_{i,d_i}
$$
the decomposition of the (non-homogeneous) equation $f_i$ into
components, homogeneous in $u_*$. (Here we somewhat abuse the
notations, using the same symbol $f_i$ both for the homogeneous
equation and for its affine presentations in the coordinates
$(z_*)$ and $(u_*)$, however this can not lead to any
misunderstanding.) The system of linear equations
$$
\Phi_{1,1}=\dots=\Phi_{k,1}=0
$$
defines the tangent space $T_pV\cong{\mathbb C}^M$.\vspace{0.1cm}

(R1) The set of all homogeneous polynomials $\Phi_{i,j}$ but the last
one $\Phi_{k,d_k}$ forms a regular sequence in
${\cal O}_{p,{\mathbb P}}$, that is, the system of homogeneous equations
$$
\{\Phi_{i,j}|_{T_pV}=0\,|\,1\leqslant i\leqslant k, 2\leqslant
j\leqslant d_i,(i,j)\neq(k,d_k)\}
$$
defines a finite set of lines in $T_pV$.\vspace{0.1cm}

{\bf Definition 0.1.} The tuple $\underline{f}\in{\cal
H}(\underline{\xi})$ is {\it regular}, if the complete
intersection $V(\underline{f})$ satisfies the conditions (R0.1,2)
at the singular point $o$ and the condition (R1) at every point
$p\neq o$.\vspace{0.1cm}

We denote the set of regular tuples by the symbol ${\cal H_{\rm
reg}}(\underline{\xi})$. The following fact is true.\vspace{0.1cm}

{\bf Theorem 0.2.} {\it ${\cal H_{\rm
reg}}(\underline{\xi})\subset{\cal H}(\underline{\xi})$ is a non-empty
Zariski open subset.}\vspace{0.1cm}

Note that the openness is obvious, so that we only need to show
that it is non-empty.\vspace{0.3cm}


{\bf 0.3. Birational rigidity and superrigidity.} The modern
method of maximal singularities goes back to the classical paper
\cite{IM}. Recall the definitions of the birational rigidity and
superrigidity for the class of varieties, considered in this paper
(for the details, see \cite[Chapter 2]{Pukh13a}). Let $X$ be a
projective rationally connected variety with ${\mathbb
Q}$-factorial terminal singularities, $D$ an effective divisor on
$X$. The {\it threshold of canonical adjunction} $c(D)=c(D,X)$ is
equal to
$$
\mathop{\rm sup}\{t\in{\mathbb Q}_+\,|\,D+tK_X\,\,\mbox{is
pseudoeffective}\}.
$$
If $\Sigma$ is a mobile linear system on $X$ and $D\in\Sigma$,
then we set $c(\Sigma)=c(D)$. The {\it virtual threshold of
canonical adjunction} $c_{\rm virt}(\Sigma,X)$ is equal to
$$
\mathop{\rm inf}_{\widetilde{X}\to
X}c(\widetilde{\Sigma},\widetilde{X}),
$$
where the infimum is taken over all birational morphisms
$\widetilde{X}\to X$ with a non-singular projective $\widetilde{X}$ and
$\widetilde{\Sigma}$ is the strict transform of the mobile system
$\Sigma$ on $\widetilde{X}$. The variety $X$ is {\it birationally
superrigid}, if $c(\Sigma)=c_{\rm virt}(\Sigma,X)$ for every
mobile linear system $\Sigma$, and {\it birationally
rigid}, if for every mobile linear system $\Sigma$ there exists a birational
self-map $\chi\in\mathop{\rm Bir}(X)$ such that the equality
$$
c(\chi^{-1}_*\Sigma)=c_{\rm virt}(\Sigma,X)
$$
holds. Now Theorem 0.1 follows from Theorem 0.2 and the following
claim.\vspace{0.1cm}

{\bf Theorem 0.3.} {\it For every tuple $\underline{f}\in{\cal
H}_{\rm reg}(\underline{\xi})$ the Fano variety
$V(\underline{f})$ is birationally superrigid.}\vspace{0.1cm}

{\bf Corollary 0.1.} {\it Let $V=V(\underline{f})$ for some
tuple $\underline{f}\in{\cal H}_{\rm
reg}(\underline{\xi})$.\vspace{0.1cm}

{\rm (i)} Assume that there is a birational map
$\chi\colon V\dashrightarrow V'$ onto the total space of the Mori
fibre space $V'\to S'$. Then $S'$ is a point and the map
$\chi$ is a biregular isomorphism.\vspace{0.1cm}

{\rm (ii)} There is no rational dominant map
$\beta\colon V\dashrightarrow S$, where $\mathop{\rm dim}S\geqslant
1$ and the general fibre is rationally connected. In particular, the
variety $V$ is non-rational.\vspace{0.1cm}

{\rm (iii)} The groups of biregular and birational automorphisms of
the variety $V$ coincide:} $\mathop{\rm Bir}V=\mathop{\rm
Aut}V$.\vspace{0.1cm}

{\bf Proof of the corollary.} These are very well known implications
of the birational superrigidity (see, for instance, \cite[Chapter
2]{Pukh13a}), taking into account that $V$ is a primitive Fano
variety.\vspace{0.3cm}


{\bf 0.4. The structure of the proof.} Theorems 0.2 and 0.3 are
independent of each other and are shown by different methods. In
Section 1 we prove Theorem 0.3. The proof is based on the
generalized $4n^2$-inequality that was recently shown in \cite{Pukh2017a}.
Assuming that the variety $V$ is not birationally superrigid, we
obtain a mobile linear system $\Sigma\subset |nH|$ on $V$
($H$ is the class of a hyperplane section of the variety
$V\subset{\mathbb P}$) with a maximal singularity $E$
(which is an exceptional divisor over $V$). The centre
$B\subset V$ of the maximal singularity $E$ is an irreducible
subvariety. If $B\neq o$ is not a singular point, then we obtain
a contradiction by word for word the same arguments as in the
non-singular case \cite{Pukh01}. If $B=o$, then by
\cite{Pukh2017a} we obtain the inequality
$$
\mathop{\rm mult}\nolimits_oZ>4n^2\mu
$$
for the self-intersection $Z=(D_1\circ D_2)$ of the linear system $\Sigma$
($D_1,D_2\in\Sigma$ are general divisors). Now the contradiction is
obtained by means of the technique of hypertangent divisors based
on the condition (R0.2). This contradiction completes the proof
of Theorem 0.3.\vspace{0.1cm}

In Section 2 we prove Theorem 0.2. The fact that a general smooth
complete intersection of the type $\underline{d}$ in ${\mathbb P}$
satisfies the condition (R1) at every point was shown in \cite{Pukh01}.
However, in our case the variety $V(\underline{f})$ with
$\underline{f}\in{\cal H}(\xi)$ has a fixed point of high
multiplicity, and for that reason we have to prove the condition (R1)
for every non-singular point $p\in V(\underline{f})$, $p\neq o$,
again. That the conditions (R0.1) and (R0.2) hold for a general tuple
$\underline{f}\in{\cal H}(\xi)$, is obvious.\vspace{0.1cm}

The proof of the condition (R1) for a non-singular point $p\in
V(\underline{f})$ is elementary but non-trivial: it requires
manipulations with coordinate presentations of the equations $f_i$
with respect to the system $z_*$ and the system $u_*$. A violation
of the regularity condition for the tuple of polynomials $\Phi_{i,j}(u_*)$
is translated into conditions for the tuple of polynomials $q_{i,j}(z_*)$
with $\xi_i\leqslant j\leqslant d_i$. In this way we prove that a general
tuple $\underline{f}\in{\cal H}(\xi)$ satisfies all regularity conditions
at all points of the variety $V(\underline{f})$.\vspace{0.3cm}


{\bf 0.5. Historical remarks and acknowledgements.} The present paper
generalizes \cite{Pukh02a}, where Fano hypersurfaces of index 1 with
a singular point of high multiplicity were shown to be birationally
rigid. Fifteen years ago the local technique of estimating the
multiplicity of the self-intersection of a mobile linear system was
weaker, and for that reason the proof given in \cite{Pukh02a} for a
point of multiplicity 3 and 4 was very hard. The generalized
$4n^2$-inequality, shown in \cite{Pukh2017a}, simplifies the
argument essentially.\vspace{0.1cm}

There are also other papers studying singular Fano varieties
from the viewpoint of their birational rigidity. After
\cite{Pukh89c} singular three-dimensional quartics became a
popular class to investigate, see \cite{CoMe,Me04,Shr08b}. Another
popular class of singular Fano varieties is formed by weighted
three-fold hypersurfaces, see \cite{CPR,Okada1,CheltsovPark2016}.
Certain families of singular higher-dimensional Fano varieties were
studied in \cite{Pukh03a,Mullany,TdF,suzuki}. The last of them is based
on the constructions of \cite{TdF}, some of which are hard to follow.
The list given above is far from being complete.\vspace{0.1cm}

There are also recent papers \cite{EP,EvP2017}, where an effective
estimate is given for the codimension of the complement to the set
of birationally superrigid varieties in the given family
(of Fano hypersurfaces with quadratic singularities, the rank of which
is bounded from below, and of Fano complete intersections of
codimension 2 with quadratic and bi-quadratic singularities, the
rank of which is bounded from below, respectively).\vspace{0.1cm}

Various technical issues, related to the constructions of this
paper, were discussed by the author in the talks given in
2012-2016 at Steklov Mathematical Institute. The author thanks the
members of divisions of Algebraic Geometry and Algebra for the
interest to his work. The author is also grateful to the
colleagues in Algebraic Geometry research group at the University
of Liverpool for the creative atmosphere and general
support.\vspace{0.1cm}

The author thanks the referee for a number of useful comments and
suggestions.


\section{Proof of birational \\ superrigidity}

In this section we prove Theorem 0.3. First (Subsection 1.1) we remind
the definition of a maximal singularity and explain, why only maximal
singularities with the centre at the singular point $o$ need to be
considered; after that (Subsection 1.2) we remind the generalized
$4n^2$-inequality that was shown in \cite{Pukh2017a}. Finally, in
Subsection 1.3 we exclude the maximal singularity.\vspace{0.3cm}

{\bf 1.1. The maximal singularity.} Fix a tuple of polynomials
$\underline{f}\in{\cal H}_{\rm reg}(\underline{\xi})$. Set
$V=V(\underline{f})$ to be the corresponding complete intersection. We have:
$$
\mathop{\rm Sing}V=\{o\},\quad\mathop{\rm Pic}V={\mathbb Z}H,\quad K_V=-H,
$$
where $H$ is the class of a hyperplane section of the variety
$V\subset{\mathbb P}$. Assume that the variety $V$ is not
birationally superrigid. This immediately implies (see
\cite[Chapter 2]{Pukh13a}) that on $V$ there is a mobile linear
system $\Sigma\subset|nH|$, $n\geqslant 1$, with a maximal
singularity: for some non-singular projective variety
$\widetilde{V}$ and a birational morphism
$\varphi\colon\widetilde{V}\to V$ there exists a
$\varphi$-exceptional prime divisor $E\subset\widetilde{V}$,
satisfying the Noether-Fano inequality
$$
\mathop{\rm ord}\nolimits_E\Sigma=
\mathop{\rm ord}\nolimits_E\varphi^*\Sigma>na(E),
$$
where $a(E)=a(E,V)\geqslant 1$ is the discrepancy of the divisor
$E$ with respect to the model $V$. Let $B=\varphi(E)\subset V$ be
the centre of the divisor $E$ on $V$. This is an irreducible
subvariety, satisfying the inequality $\mathop{\rm
mult}_B\Sigma>n$. By the Lefschetz theorem for the numerical Chow
group of algebraic cycles of codimension 2 on $V$ we have the
equality $A^2V={\mathbb Z}H^2$. Now we are able to exclude the
simplest case $\mathop{\rm codim}(B\subset V)=2$.\vspace{0.1cm}

Indeed, if $\mathop{\rm codim}(B\subset V)=2$, then $B\sim mH^2$
for some $m\geqslant 1$. Consider the self-intersection
$Z=(D_1\circ D_2)$ of the linear system $\Sigma$. Obviously, $Z\sim n^2H^2$.
On the other hand, $Z=\beta B+Z_1$, where $\beta>n^2$ and $Z_1$ is an
effective cycle of codimension 2 that does not contain $B$ as a component.
Taking the classes of the cycles in $A^2V$, we obtain the inequality
$n^2\geqslant\beta m>mn^2$. The contradiction excludes the case
$\mathop{\rm codim}(B\subset V)=2$.\vspace{0.1cm}

If $\mathop{\rm codim}(B\subset V)\geqslant 3$ and $B\neq o$, then
the inequality $\mathop{\rm mult}_BZ>4n^2$ holds (the classical
$4n^2$-inequality that goes back to \cite{IM}, see \cite[Chapter 2]{Pukh13a}).
Repeating the arguments of \cite[Section 2]{Pukh01} word for word, we get
a contradiction, excluding the case under consideration. We are able to repeat
the arguments word for word, because the proof in
\cite[Section 2]{Pukh01}, is based on the regularity condition, which is
identical to the condition (R1).\vspace{0.1cm}

So we are left with the only option: $B=o$. Excluding this option,
we complete the proof of Theorem 0.3.\vspace{0.3cm}


{\bf 1.2. The generalized $4n^2$-inequality.} Recall that
$\underline{\mu}=(\mu_1,\dots,\mu_l)$ with
$2\leqslant\mu_1\leqslant\dots\leqslant\mu_l$ is the type
of the complete intersection singularity $o\in V$, so that by the
condition (R0.1)
$$
\mu=\mu_1\cdots\mu_l=\mathop{\rm mult}\nolimits_oV.
$$
Again, let $Z=(D_1\circ D_2)$ be the self-intersection of the linear
system $\Sigma$. The following fact is true.\vspace{0.1cm}

{\bf Theorem 1.1.} {\it The following inequality holds:}
$$
\mathop{\rm mult}\nolimits_oZ>4n^2\mu.
$$

{\bf Proof}: this is a particular case of the theorem shown in
\cite{Pukh2017a} under essentially weaker assumptions about the
generality of the singularity $o\in V$. Q.E.D.\vspace{0.3cm}


{\bf 1.3. Hypertangent divisors.} Now let us use the technique of
hypertangent linear systems. We use the notations of Subsection
0.2 and work in the affine chart ${\mathbb C}^{M+k}$ of the space
${\mathbb P}$ with the coordinates $z_1,\dots,z_{M+k}$. Let
$j\geqslant 1$ be an integer.\vspace{0.1cm}

{\bf Definition 1.1.} The linear system
$$
\Lambda(j)=
\left\{\left.\left(\sum^k_{i=1}\sum^{d_i-1}_{\alpha=\xi_i}
f_{i,\alpha}s_{i,j-\alpha}\right)\right|_V =0\right\},
$$
where $s_{i,j-\alpha}$ independently run through the set
of homogeneous polynomials of degree $j-\alpha$ in the variables
$z_*$ (if $j-\alpha<0$, then $s_{j-\alpha}=0$), and
$$
f_{i,\alpha}=q_{i,\xi_i}+\dots+q_{i,\alpha}
$$
is the truncated equation $f_i$, is called the {\it $j$-th
hypertangent system at the point} $o$.\vspace{0.1cm}

For $j\geqslant 1$ set
$$
c(j)=\sharp\{(i,\alpha)\,|\,\, i=1,\dots,k,
\xi_i\leqslant\alpha\leqslant \min\{j,d_i-1\}\,\}.
$$
Set also $a=\mathop{\rm min}\{j\geqslant 1\,|\,\,c(j)\geqslant 1\}$.
Obviously, the system $\Lambda(j)$ is non-empty only for $j\geqslant a$.
Moreover, the equality
$$
\mathop{\rm codim}\nolimits_o(\mathop{\rm Bs}\Lambda(j)\subset
V)=c(j)
$$
holds. Now set $m(j)=c(j)-c(j-1)$ and choose in every non-empty
hypertangent system $\Lambda(j)$ precisely $m(j)$ general divisors
$$
D_{j,1},\dots,D_{j,m(j)}
$$
(if $m(j)=0$, then we do not choose any divisors in the system
$\Lambda(j)$), for $j=a,\dots,d_k-1$. It is easy to see that
we obtained a tuple consisting of
$$
m=\sum^k_{i=1}(d_i-\xi_i)
$$
effective divisors on $V$. Denoting by the symbol
$|D_{j,\alpha}|$ the support of the divisor $D_{j,\alpha}$,
we conclude: by the condition (R0.2)
$$
\mathop{\rm codim}\nolimits_o
\left(\left(\bigcap^{d_k-1}_{j=a}\bigcap^{m(j)}_{\alpha=1}
|D_{j,\alpha}|\right)\subset V\right)=m
$$
(the symbol $\mathop{\rm codim}_o$ stands for the codimension in a
neighborhood of the point $o$). Let us place the divisors $D_{j,\alpha}$ in
the {\it standard order}: $D_{j_1,\alpha_1}$ precedes
$D_{j_2,\alpha_2}$, if $j_1<j_2$ or $j_1=j_2$, but
$\alpha_1<\alpha_2$. We obtain a sequence
$$
R_1,\dots,R_m
$$
of effective divisors on $V$. Now let us consider the effective cycle
of the self-intersection $Z=(D_1\circ D_2)$ of the mobile system
$\Sigma$. We denote its support by the symbol $|Z|$. By the
condition (R0.2) for every $i\in\{3,\dots,m\}$ we have the equality
$$
\mathop{\rm codim}\nolimits_o\left(\bigcap^i_{j=3}|R_j|\cap|Z|\right)=i.
$$

Let $Y\subset|Z|$ be an irreducible component of the support $|Z|$,
which has the maximal value of the ratio
$\mathop{\rm mult}_o/\mathop{\rm deg}$ of the multiplicity at the point $o$
to the degree (in the space ${\mathbb P}$). By Theorem 1.1
$$
\frac{\mathop{\rm mult}_o}{\mathop{\rm deg}}Y>\frac{4n^2}{n^2d}\mu=\frac{4}{d}\mu,
$$
where $d=d_1\cdots d_k=\mathop{\rm deg}V$ and
$\mu=\mu_1\cdots\mu_l=\xi_1\cdots\xi_k$. Now let us construct in
the usual way (see \cite[Chapter 3]{Pukh13a}) a sequence of irreducible
subvarieties $Y_2=Y,Y_3,\dots,Y_m$, satisfying the following properties:
\begin{itemize}

\item $\mathop{\rm codim}(Y_i\subset V)=i$,

\item $Y_{i+1}$ is an irreducible component of the algebraic cycle
of the scheme-theoretic intersection $(Y_i\circ R_{i+1})$ with the
maximal value of $\mathop{\rm mult}_o/\mathop{\rm deg}$.
\end{itemize}

Such construction is possible, because by the condition (R0.2) and
the genericity of the divisors $D_{j,\alpha}$ in the linear system
$\Lambda(j)$ we have
$$
Y_i\not\subset|R_{i+1}|.
$$
The irreducible subvariety $Y_m$ is of positive dimension and
contains the point $o$. Let us see how the value of the ratio
$\mathop{\rm mult}_o/\mathop{\rm deg}$ changes when we make the step
from $Y_i$ to $Y_{i+1}$.\vspace{0.1cm}

Let $R_i=D_{b,\alpha}$ for some $b\in\{a,\dots,d_k-1\}$ and
$\alpha\in\{1,\dots,m(b)\}$, in particular, $R_i\in\Lambda(b)$.
The number
$$
\beta_i=\beta(R_i)=\frac{b+1}{b}
$$
is called the {\it slope} of the divisor $R_i$.\vspace{0.1cm}

{\bf Proposition 1.1.} {\it For $i=3,\dots,m$ the following inequality holds:}
$$
\frac{\mathop{\rm mult}_o}{\mathop{\rm deg}}Y_i\geqslant\beta_i\,\,
\frac{{\mathop{\rm mult}_o}}{{\mathop{\rm deg}}}Y_{i-1}.
$$

{\bf Proof.} By construction, $D_{b,\alpha}\sim bH$. For the
strict transform $D^+_{b,\alpha}$ on $V^+$ we have
$D^+_{b,\alpha}\sim bH-(b+1)Q$, since
$$
f_{j,b}|_V=(q_{j,\xi_j}+\dots+q_{j,b})|_V=(-q_{j,b+1}+\dots)|_V
$$
for every $j\in\{1,\dots,k\}$ such that
$\xi_j\leqslant b\leqslant d_j-1$. Q.E.D.
for the proposition.\vspace{0.1cm}

The maximal possible value of the slope $\beta_i$ is 2
(for $b=1$). The hypertangent divisors are ordered in such way that
their slopes do not increase:
$\beta_1\geqslant\beta_2\geqslant\dots\geqslant\beta_m$.
Therefore,
$$
\frac{\mathop{\rm mult}_o}{\mathop{\rm deg}}Y_m>
4\frac{\xi_1\dots\xi_k}{d_1\dots d_k}
\prod^m_{i=3}\beta_i\geqslant
\frac{4}{\beta_1\beta_2}\frac{\xi_1\dots\xi_k}{d_1\dots d_k}
\prod^m_{i=1}\beta_i\geqslant
\frac{\xi_1\dots\xi_k}{d_1\dots d_k}
\prod^k_{i=1}\prod^{d_i-1}_{j=\xi_i}\frac{j+1}{j}=1.
$$
Thus $\mathop{\rm mult}_oY_m>\mathop{\rm deg}Y_m$, which is, of course,
impossible. We obtained a contradiction which completes the proof
of Theorem 0.3.


\section{Regular complete intersections}

In this section we prove Theorem 0.2. In Subsection 2.1 we
consider a convenient system of coordinates for the points $o$ and
$p\neq o$ and give the formulas that relate the local
presentations of the polynomials $f_i$ at those points. In
Subsection 2.2 we explain what impact the position of the line
connecting the points $o$ and $p$ has for the regularity
conditions. In Subsection 2.3 we reduce the claim of Theorem 0.2
to a certain fact about tuples of polynomials $q_{i,j}$ (the
homogeneous components of the equations $f_i$). In Subsection 2.4
we consider the special case of violation of the regularity
conditions when the polynomials $\Phi_{i,j}$ vanish on a line.
Finally, in Subsection 2.5 we consider the general case of
violation of the regularity conditions and complete the proof of
Theorem 0.2.\vspace{0.3cm}

{\bf 2.1. Preliminary constructions.} We fix the system of affine coordinates
$(z_*)=(z_1,\dots, z_{M+k})$ with the origin at the point $o$. The
symbol ${\cal P}_{l,M+k}$ stands for the linear space of homogeneous
polynomials of degree $l$ in $M+k$ variables $(z_*)$. Set
$$
{\cal P}_{[a,b],M+k}=\prod^b_{l=a}{\cal P}_{l,M+k}
$$
to be the space of polynomials of the form $q_a+\dots +q_b$, where
$q_l\in {\cal P}_{l,M+k}$. Finally, let
$$
{\cal P}=\prod^k_{i=1}{\cal P}_{[\xi_i,d_i],M+k}
$$
be the space of tuples $\underline{f}=(f_1,\dots ,f_k)$. Obviously,
for a general tuple $\underline{f}\in {\cal P}$ the conditions (R0.1) and
(R0.2) are satisfied.\vspace{0.1cm}

Let $p\in {\mathbb P}$, $p\neq o$, be an arbitrary point. We will assume
that $p$ lies in the affine chart with coordinates $(z_*)$ and, moreover, has in
that chart the coordinates $(1,0,\dots,0)$. Let
$$
u_1=z_1-1, u_2=z_2,\dots,u_{M+k}=z_{M+k}
$$
be a system of coordinates with the origin at the point $p$.
It can be written as
$u_*=(u_1,\dots,u_{M+k})=(u_1,z_2,\dots,z_{M+k})$. Recall that the equation
$f_i$ can be naturally written in the form
$$
f_i=q_{i,\xi_i}+\dots+q_{i,d_i}.
$$
Set $q_{i,j}=q_{i,j,j}+z_1q_{i,j,j-1}+\dots+z^j_1q_{i,j,0}$, where
$q_{i,j,l}$ is a homogeneous polynomial of degree $l\leqslant j$
in the variables $z_2,\dots,z_{M+k}$. In the new coordinates
$(u_*)$ the polynomial $f_i$ takes the form
\begin{equation}\label{17.03.2017.1}
f_i=\sum^{d_i}_{e=0}\left[\sum^e_{\alpha=0}u^{e-\alpha}_1
\left(\sum^{d_i}_{j=\mathop{\rm max}(\xi_i,e)} {j-\alpha\choose
e-\alpha}q_{i,j,\alpha}\right)\right].
\end{equation}
The expression in the square brackets is the homogeneous polynomial
$\Phi_{i,e}$ in the notations of Subsection 0.2. Therefore, the shift of the first
coordinate $z_1=1+u_1$, with the other coordinates being the
same, defines a linear map
$$
\tau\colon{\cal P}\to\prod^k_{i=1}{\cal P}_{[0,d_i],M+k}(u_*).
$$
The following claim is true.\vspace{0.1cm}

{\bf Theorem 2.1.} {\it The set of tuples $\underline{f}\in{\cal P}$,
such that
$$
\tau(\underline{f})\in\prod^k_{i=1}{\cal P}_{[1,d_i],M+k}(u_*)
$$
(that is, $p\in V(\underline{f})$) and $\tau(\underline{f}$) does not
satisfy the condition (R1), is of codimension at least
$M+k+1$ in the space} ${\cal P}$.\vspace{0.1cm}

{\bf Proof of Theorem 0.2.} Theorem 2.1 implies that the set of tuples
$\underline{f}\in{\cal H}(\underline{\xi})$ (here we again abuse the
notations, using the same symbol $f_i$ both for a homogeneous polynomial
of degree $d_i$ on ${\mathbb P}$ and for its non-homogeneous presentation
in the coordinates $(z_*)$) such that the variety $V(f)$ is not regular
at at least one point $p\in V(\underline{f})$, $p\neq o$, has a positive
codimension in ${\cal H}(\underline{\xi})$. Therefore, the set
${\cal H}_{\rm reg}(\underline{\xi})$ is non-empty. It is obviously
open. Q.E.D. for Theorem 0.2.\vspace{0.1cm}

The rest of this section is a proof of Theorem 2.1.\vspace{0.1cm}

Note that $\Phi_{i,0}=q_{i,\xi_i,0}+\dots +q_{i,d_i,0}$ and the
equalities $\Phi_{i,0}=0$, $i=1,\dots, k$, give $k$ independent
linear conditions for the polynomials $f_i$ (expressing the fact that
$p\in V(\underline{f})$). Furthermore, for the linear terms we get
$$
\Phi_{i,1}=(\xi_i q_{i,\xi_i,0}+\dots +d_i q_{i,d_i,0})u_1+q_{i,\xi_i,1}+
\dots +q_{i,d_i,1}
$$
and the linear dependence of these linear forms gives in addition
$M+1$ independent conditions for the coefficients of the
polynomials $f_i$. It is for this reason that the open set ${\cal
H}(\underline{\xi})$ is non-empty.\vspace{0.1cm}

For $k$ fixed linearly independent linear forms $L_1,\dots ,L_k$
in the variables $u_*$ we denote by the symbol
$$
{\cal P}[p;L_1,\dots ,L_k]\subset {\cal P}
$$
the affine subspace of tuples $\underline{f}$ such that
$\Phi_{1,0}=\dots =\Phi_{k,0}=0$ and
$\Phi_{1,1}=L_1$,\dots , $\Phi_{k,1}=L_k$.\vspace{0.3cm}


{\bf 2.2. The line connecting the points $o$ and $p$.} Let us
denote this line by the symbol $[o,p]$. We say that we are in the
{\it non-special case}, if $[o,p]\not\subset T_pV$, and in the
{\it special case}, if $[o,p]\subset T_pV$. Consider first the
non-special case.\vspace{0.1cm}

In the coordinates $u_*$ on the space
$T_p{\mathbb P}\cong{\mathbb C}^{M+k}$ the line
$[o,p]$ is given by the equations
$$
u_2=\dots=u_{M+k}=0.
$$
If $[o,p]\not\subset T_pV$, then
$$
\mathop{\rm dim}\langle u_2|_{T_pV},\dots,u_{M+k}|_{T_pV}\rangle=M,
$$
so that for every $j\geqslant 1$ the space
$$
\{q_j|_{T_pV}\,\,|\,\, q_j\in{\cal P}_{j,M+k-1}\}
$$
(where ${\cal P}_{j,M+k-1}$ is the linear space of
homogeneous polynomials of degree $j$ in $u_2,\dots,u_{M+k}$)
is the whole space of homogeneous polynomials of degree $j$
on $T_pV$. It follows from (\ref{17.03.2017.1}) that
$\Phi_{i,e}=q_{i,d_i,e}+(*)$, where $(*)$ is a
linear combination of terms
$u_1^{e-\alpha}q_{i,j,\alpha}$ with either
$j<d_i$ or $j=d_i$ and $\alpha<e$. As the polynomials
$q_{i,j,\alpha}$ are arbitrary homogeneous polynomials
of degree $\alpha$ in $u_2,\dots,u_{M+k}$, we conclude,
that when $\underline{f}$ runs through the space
${\cal P}[p; L_1,\dots,L_k]$, the homogeneous polynomials
$\Phi_{i,e}|_{T_pV},e\geqslant 2$, run through the whole
spaces of homogeneous polynomials of degree $e$ on
$T_pV$, independently of each other. For that reason, by
\cite[Section 3]{Pukh01} the set of tuples
$\underline{f}\in{\cal P}[p;L_*]$, for which the condition
(R1) is violated at the point $p$, has in ${\cal P}[p,L_*]$
the codimension at least $M+1$. Since for different tuples
$(L_1,\dots,L_k)\neq(L'_1,\dots,L'_k)$
of linearly independent linear forms the affine subspaces
${\cal P}[p;L_*]$ and ${\cal P}[p;L'_*]$ are disjoint, the reference
to \cite[Section 3]{Pukh01} completes the proof of Theorem 2.1
in the non-special case.\vspace{0.1cm}

Starting from this moment, we assume that we are in the
special case, that is, for all $i=1,\dots,k$ we have
$L_i|_{[o,p]}\equiv 0$. Explicitly, this means that
for all $i=1,\dots,k$ the equality
$$
\xi_iq_{i,\xi_i,0}+\dots+d_iq_{i,d_i,0}=0
$$
holds. If $\xi_i\leqslant d_i-1$, that is, the hypersurface
$\{f_i=0\}$ is not a cone with the vertex at the point $o$,
we obtain a new independent condition for $\underline{f}$.
Therefore, the set
$$
\bigsqcup_{\rm special}{\cal P}[p; L_1,\dots,L_k],
$$
where the union is taken over all tuples $(L_*)$, consisting of
$k$ linearly independent forms, vanishing on the line $[o,p]$,
has in ${\cal P}$ the codimension
$$
k+\sharp\{i\,|\,\xi_i\leqslant d_i-1,\,\, i=1,\dots,k\}.
$$

Set ${\mathbb T}={\mathbb P}(T_pV)\cong {\mathbb P}^{M-1}$.
The point, corresponding to the line $[o,p]$, we denote by the symbol
$\omega$. Theorem 2.1 is implied by the following claim.\vspace{0.1cm}

{\bf Proposition 2.1.} {\it In the special case the set of
tuples $\underline{f}\in {\cal P}[p; L_*]$ such that the system of
equations
\begin{equation}\label{18.03.2017.1}
\Phi_{i,j}|_{\mathbb T}=0,\quad 1\leqslant i\leqslant k,\,\,
2\leqslant j\leqslant d_i,\,\, (i,j)\neq (k,d_k)
\end{equation}
has in ${\mathbb T}$ a positive-dimensional set of solutions,
is of codimension at least
$$
M+1-\sharp\{i\,|\,\xi_i\leqslant d_i-1,\,\, i=1,\dots,k\}
$$
in the space} ${\cal P}[p; L_*]$.\vspace{0.3cm}


{\bf 2.3. Plan of the proof of Proposition 2.1.}
The linear forms $\Phi_{i,1}=L_i$ and the projective space
${\mathbb T}$ are fixed. Placing the polynomials
$\Phi_{i,j}|_{\mathbb T}$ in the standard order
(which means that $(i_1,j_1)<(i_2,j_2)$, if $i_1<i_2$ or $i_1=i_2$, but
$j_1<j_2$), we get $M-1$ polynomials on
${\mathbb P}^{M-1}$:
$$
p_1,p_2,\dots,p_{M-1},
$$
$\mathop{\rm deg}p_{i+1}\geqslant\mathop{\rm deg}p_i$. In the
special case it is not true that $p_i$ run through the
corresponding spaces of polynomials independently of each other,
and for that reason the estimates that were obtained in
\cite[Section 3]{Pukh01} can not be applied directly. Let us
consider the following subsets of the affine space ${\cal A}={\cal
P}[p; L_*]$.\vspace{0.1cm}

Let ${\cal B}_{\rm line}\subset{\cal A}$ be the set of tuples
$\underline{f}\in{\cal A}$ such that
$p_i|_R\equiv 0$ for some line $R\subset{\mathbb T}$ for
all $i=1,\dots,M-1$. Furthermore, set
${\cal B}_i\subset{\cal A}\backslash{\cal B}_{\rm line}$,
where $i=1,\dots,M-1$, to be the set of tuples
$\underline{f}\in{\cal A}$ such that
$$
\mathop{\rm codim}(\{p_1=\dots=p_{i-1}=0\subset{\mathbb T}\})=i-1,
$$
but for some irreducible component $B$ of the set
$\{p_1=\dots=p_{i-1}=0\}$ we have $p_i|_B\equiv 0$. (For $i=1$
this condition means that $p_1\equiv 0$.) Recall that
$$
c_*=\sharp\{i\,\,|\,\,\xi_i=d_i,\,\,i=1,\dots,k\}.
$$

Proposition 2.1 is implied by the following two facts.\vspace{0.1cm}

{\bf Proposition 2.2.} {\it The following inequality holds:}
$$
\mathop{\rm codim}({\cal B}_{\rm line}\subset{\cal A})\geqslant
M+1+c_*-k.
$$

{\bf Proposition 2.3.}  {\it For all $i=1,\dots,M-1$ the following inequality holds:}
$$
\mathop{\rm codim}({\cal B}_i\subset{\cal A})\geqslant M+1.
$$

{\bf Remark 2.1.} Let $(v_0:v_1:\dots : v_{M-1})$ be a system of
homogeneous coordinates on ${\mathbb T}$, in which the point
$\omega$ is the point $(1:0:\dots:0)$. The formula
(\ref{17.03.2017.1}) implies that for fixed polynomials
$p_1,\dots,p_{i-1}$ the set, which the polynomial $p_i$ runs
through, is a disjoint union of affine subspaces of the form
$$
p^{\circ}_i+{\cal P}_{\mathop{\rm deg} p_i,M-1}(v_1,\dots,v_{M-1}),
$$
where $p^{\circ}_i$ is some polynomial. Applying the method
of linear projections (see \cite[Chapter 3, Section 1]{Pukh13a}),
we obtain the inequality
$$
\mathop{\rm codim}({\cal B}_i\subset{\cal A}) \geqslant
{M-i-1+\mathop{\rm deg}p_i\choose\mathop{\rm deg}p_i}.
$$
For high values of $i$ (those close to $M-1$) this estimate can be
not strong enough. However, for $i=1,\dots,k$ it gives more than
we need:
$$
\mathop{\rm codim}({\cal B}_i\subset{\cal A})\geqslant
{M-k+1\choose 2}.
$$
Therefore it is sufficient to prove Proposition 2.3 for
$i\geqslant k+1$, so that $\mathop{\rm deg}p_i \geqslant 3$.\vspace{0.1cm}

The proof of Proposition 2.2 is given in Subsection 2.4, of Proposition 2.3 in
Subsection 2.5.\vspace{0.3cm}


{\bf 2.4. The case of a line.} First of all, let us break the set
${\cal B}_{\rm line}$ into two (overlapping) subsets, corresponding
to the case when the line $R$ contains the point
$\omega$ (the subset ${\cal B}^+_{\rm line}$) and when it does not
(the subset ${\cal B}^-_{\rm line}$). Let us estimate the
codimensions of these sets in the space ${\cal A}$ separately.\vspace{0.1cm}

For $\underline{f}\in{\cal B}^-_{\rm line}$ the conditions
$p_i|_R\equiv 0$, $i=1,\dots,M-1$, similarly to the non-special
case, give $\sum\limits^{M-1}_{i=1}(\mathop{\rm deg}p_i+1)$
independent conditions for $\underline{f}$. Since the line $R$
varies in a $2(M-2)$-dimensional family, we obtain the estimate
$$
\mathop{\rm codim}({\cal B}^-_{\rm line}\subset{\cal A})\geqslant
\sum^{M-1}_{i=1}\mathop{\rm deg}p_i-M+3.
$$

{\bf Lemma 2.1.} {\it The following inequality holds:}
$$
\sum^{M-1}_{i=1}\mathop{\rm deg}p_i\geqslant 2M-2.
$$

{\bf Proof.} It is easy to check that
\begin{equation}\label{21.03.2017.1}
\sum^{M-1}_{i=1}\mathop{\rm deg}p_i=\sum^{k-1}_{i=1}\frac{d_i(d_i+1)}{2}+
\frac{(d_k-1)d_k}{2}-k.
\end{equation}
Furthermore, it is easy to see that if
$d_{i-1}\leqslant d_i-2$, then, replacing $d_{i-1}$ by
$d_{i-1}+1$, and $d_i$ by  $d_i-1$, we will not increase the
value of the expression (\ref{21.03.2017.1}). Therefore, the minimum
of that expression for the fixed $k$ and $M$ is attained at
$$
d_1=\dots=d_{k-l}=a+1,\,\,d_{k-l+1}=\dots=d_k=a+2,
$$
where $M=ka+l$, $l\in\{0,1,\dots,k-1\}$. Now elementary
computations complete the proof of the lemma. Q.E.D.\vspace{0.1cm}

Lemma 2.1 implies the estimate
$$
\mathop{\rm codim}({\cal B}^-_{\rm line}\subset{\cal A})\geqslant M+1,
$$
which is stronger than the inequality of Proposition 2.2. For that reason,
Proposition 2.2 is implied by the following claim.\vspace{0.1cm}

{\bf Proposition 2.4.}  {\it The following inequality is true:}
$$
\mathop{\rm codim}({\cal B}^+_{\rm line}\subset{\cal A})\geqslant M+1+c_*-k.
$$

{\bf Proof.} Let $R\ni\omega$ be a line. In the notations of Remark 2.1 let
$$
\lambda=(0:a_1:\dots:a_{M-1})=R\cap\{v_0=0\}.
$$
The conditions $p_i|_R\equiv 0$, $i=1,\dots,M-1$, give a smaller
codimension than in the case $R\not\ni\omega$, considered above.
As a compensation, the lines $R\ni\omega$ vary in a
$(M-2)$-dimensional family. Let us fix the line $R$ and the point
$\lambda$.\vspace{0.1cm}

{\bf Lemma 2.2.} {\it For every $i=1,\dots,k-1$ the conditions
$$
\Phi_{i,2}|_R\equiv\dots\equiv\Phi_{i,d_i}|_R\equiv 0
$$
are equivalent to the conditions}
$$
q_{i,j,\alpha}(\lambda)=0,
$$
$\xi_i\leqslant j\leqslant d_i$,
$\alpha=0,1,\dots,j$.\vspace{0.1cm}

{\bf Proof.} For the homogeneous polynomial
$$
\Phi(v_*)=v^l_0q_0+v^{l-1}_0q_1+\dots+v_0q_{l-1}+q_l,
$$
where $q_i(v_1,\dots,v_{M-1})$ is a homogeneous polynomial
of degree $i$, the condition $\Phi|_R\equiv 0$ means that
$$
q_0=q_1(\lambda)=\dots=q_l(\lambda)=0.
$$
The formula (\ref{17.03.2017.1}) implies that if all polynomials
$\Phi_{i,j}$ (for a fixed $i\leqslant k-1$) vanish identically on the
line $R$, then the equalities
$$
\sum^{d_i}_{j=\mathop{\rm max}(\xi_i,e)}{j-\alpha\choose
e-\alpha}q_{i,j,\alpha}(\lambda)=0
$$
hold for all $e=0,\dots,d_i$ and $\alpha=0,\dots,e$. Setting
$e=d_i$, we obtan the system of equalities
$$
q_{i,d_i,\alpha}(\lambda)=0,\,\,\alpha=0,\dots,d_i.
$$
If $\xi_i=d_i$, then the claim of the lemma is shown.\vspace{0.1cm}

Assume that $\xi_i\leqslant d_i-1$. Setting
$e=d_i-1$, we obtain the system of equalities
$$
q_{i,d_i-1,\alpha}(\lambda)+{d_i-\alpha\choose d_i-1-\alpha}q_{i,d_i,\alpha}(\lambda)=0
$$
for $\alpha=0,\dots,d_i-1$, whence, taking into account the
previous equalities, we conclude that
$$
q_{i,d_i-1,\alpha}(\lambda)=0,\,\,\alpha=0,\dots,d_i-1.
$$
Proceeding in the same spirit, we consider the values
$e=d_i-2,\dots,\xi_i$ and complete the proof of the lemma.
Q.E.D.\vspace{0.1cm}

{\bf Lemma 2.3.} {\it The conditions
$$
\Phi_{k,2}|_R\equiv\dots\equiv\Phi_{k,d_k-1}|_R\equiv 0
$$
define a linear subspace of codimension
$$
\frac12[d_k(d_k+1)-\xi_k(\xi_k+1)]
$$
in the space of tuples of homogeneous polynomials}
$q_{k,j,\alpha},\,\,\xi_k\leqslant j\leqslant d_k-1,\,\,\alpha=0,1,\dots,j$.\vspace{0.1cm}

{\bf Proof.} Adding the condition
$$
\Phi_{k,d_k}|_R\equiv 0
$$
and applying the previous lemma, we obtain
$\frac12[(d_k+1)(d_k+2)-\xi_k(\xi_k+1)]$ independent linear
conditions $q_{k,j,\alpha}(\lambda)=0$. Vanishing of the form
$\Phi_{k,d_k}$ on the line $R$ adds precisely $d_k+1$ linear conditions.
Q.E.D. for the lemma.\vspace{0.1cm}

Combining Lemmas 2.2 and 2.3, we see that Proposition
2.4 follows from the inequality
$$
\frac12\sum^k_{i=1}[({d_i+1})(d_i+2)-\xi_i(\xi_i+1)]-(d_k+1)-(M-2)\geqslant
M+1+c_*-k.
$$
The last inequality is precisely (\ref{07.03.2017.1}). Q.E.D. for
Propositions 2.4 and 2.2.\vspace{0.3cm}


{\bf 2.5. Estimating the codimension in the general case.} Let us
show Proposition 2.3. By Remark 2.1 we assume that $i\geqslant
k+1$, so that $\mathop{\rm deg}p_i\geqslant 3$. We use the method
of good sequences and associated subvarieties, developed in
\cite[Section 3]{Pukh01}, see also \cite[Chapter 3, Section
3]{Pukh13a}.\vspace{0.1cm}

Let ${\cal B}_{i,b}\subset{\cal B}_i$ be the subset of tuples
$\underline{f}\in{\cal A}$ such that for some irreducible
component $B$ of the set $\{p_1=\dots=p_{i-1}=0\}$ (which has
codimension $i-1$ in ${\mathbb T}$, since $\underline{f}\in{\cal
B}_i$), such that $\mathop{\rm codim}(\langle
B\rangle\subset{\mathbb T})=b$, we have $p_i|_B\equiv 0$. The
parameter $b$ runs through the set of values $\{0,1,\dots,i-1\}$
for $i\leqslant M-2$, and through the set $\{0,\dots,M-3\}$ for
$i=M-1$. When $b=i-1$, the component $B$ is a linear subspace in
${\mathbb T}$, and the codimension $\mathop{\rm codim}({\cal
B}_{i,i-1}\subset{\cal A})$ can be calculated explicitly in the
same way as the codimension of the subset ${\cal B}_{\rm line}$,
but we do not need that.\vspace{0.1cm}

In order to estimate the codimension of the subsets ${\cal
B}_{i,b}$, we need a small modificati\-on of the technique of
\cite[Section 3]{Pukh01}. Let $P$ be a linear subspace of
codimension $b$ in ${\mathbb T}$. By the symbol ${\cal
B}_{i,b}(P)$ we denote the subset of tuples $\underline{f}\in
{\cal B}_{i,b}$ such that the linear span of the irreducible
subvariety $B$ (see the definition of the set ${\cal B}_{i,b}$) is
$P$. Obviously,
$$
\mathop{\rm codim}({\cal B}_{i,b}\subset{\cal A})\geqslant
\mathop{\rm codim}({\cal B}_{i,b}(P)\subset{\cal A})-b(M-b).
$$
Furthermore, for a subset of indices
$$
I=\{j_1<\dots<j_{i-1-b}\}\subset\{1,\dots,i-1\}
$$
let ${\cal B}_{i,b,I}(P)\subset{\cal B}_{i,b}(P)$ be the subset of
tuples $\underline{f}\in{\cal B}_{i,b}(P)$ such that there exists
a sequence of irreducible subvarieties
$$
Y_0=P,\,Y_1,\,\dots,\,Y_{i-1-b}=B,
$$
satisfying the following properties:
\begin{itemize}

\item for every $l\in\{1,\dots,i-1-b\}$ and every index
$j_{l-1}<j<j_l$ (where $j_0=0$) the polynomial $p_j$ vanishes
identically on $Y_{l-1}$,

\item for every $l\in\{1,\dots,i-1-b\}$ we have
$p_{j_l}|_{Y_{l-1}}\not\equiv 0$ and $Y_l\subset Y_{l-1}$
is an irreducible component of the closed set
$\{p_{j_l}|_{Y_{l-1}}=0\}$, containing the subvariety $B$.
\end{itemize}
In the terminology of \cite[Section 3]{Pukh01} the polynomials
$p_{j_l}|_P$, $l=1,\dots,i-1-b$, form a good sequence and $B$ is
one of its associated subvarieties. Obviously,
$$
{\cal B}_{i,b}(P)=\bigcup_I{\cal B}_{i,b,I}(P).
$$

{\bf Lemma 2.4.} {\it The following inequality holds:}
$$
\mathop{\rm codim}({\cal B}_{i,b,I}(P)\subset{\cal A})
\geqslant(2b+3)(M-1-b)-2.
$$

{\bf Proof} repeats the arguments in \cite[Section 3]{Pukh01} (the
proof of Proposition 4) or in \cite[Chapter 3, Section
3]{Pukh13a}, but $(M-1)$ needs to be replaced by $(M-2)$, since as
we have already mentioned, $\Phi_{i,e}=q_{i,d_i,e}+\dots$, where
$q_{i,d_i,e}$ is an arbitrary homogeneous polynomial of degree $e$
in the variables $u_2,\dots,u_{M+k}$. We check the polynomials
$p_j$ that were not included into the good sequence, one by one.
When the polynomials $p_{\gamma}$ with $\gamma<j$ are fixed, the
condition $p_j|_{Y_{l-1}}\equiv 0$ imposes on the coefficients of
the polynomial $p_j$ at least
$$
\mathop{\rm deg} p_j (M-2-b)+1\geqslant 2(M-2-b)+1
$$
independent conditions, since $\langle Y_{l-1}\rangle=P$
(recall that $Y_{l-1}\supset B$). The condition $p_i|_B\equiv 0$
gives (with $p_1,\dots,p_{i-1}$ fixed) at least
$$
\mathop{\rm deg}p_i(M-2-b)+1\geqslant 3(M-2-b)+1
$$
independent conditions. Putting together, we complete
the proof of the lemma. Q.E.D.\vspace{0.1cm}

Now we complete the proof of Proposition 2.3. Let us look first at
the values $i\leqslant M-2$, when the parameter $b$ takes the
values $0,1,\dots, i-1$. Let us consider the quadratic function
$$
\varphi_1(t)=(2t+3)(M-1-t)-t(M-t).
$$
Since $\varphi''_1(t)=-1<0$, its minimum on the set
$[0,i-1]$ is attained either at $t=0$, or at $t=i-1$.
Therefore, for $i=k+1,\dots,M-2$ we get
$$
\mathop{\rm codim}({\cal B}_i\subset{\cal A})
\geqslant\mathop{\rm min}\{3M-5,(M-i-1)(i+2)+1\}.
$$
Since $3M-5\geqslant M+1$, which is what we need, let us consider
the quadratic function
$$
\varphi_2(t)=(M-t-1)(t+2)+1.
$$
Again, $\varphi''_2(t)=-1<0$, so that its minimum on the set
$[k+1,M-2]$ is attained either at $t=k+1$, or at $t=M-2$. In the
latter case we get $\varphi_2(M-2)=M+1$, as required. Therefore,
in order to prove Proposition 2.3 for $i\leqslant M-2$, it is
sufficient to show the inequality
$$
(k+3)(M-k-2)+1\geqslant M+1.
$$
For $M$ fixed, the left hand side of the last inequality is a
quadratic function in $k$ with $-k^2$ as the senior term.
Therefore, its minimum is attained either when $k=2$, when the
value of the left hand side is $5M-19$ (which is not less than
$M+1$), or for the maximal possible value of $k$ (for the given
$M$). If $M$ is odd, then we need to check the value $k=\frac12
(M-1)$ or, equivalently, substitute $M=2k+1$ into both left and
right hand sides of the last displayed inequality which gives
$$
(k+3)(k-1)+1\geqslant 2k+2,
$$
which is true for $k\geqslant 2$. If $M$ is even, we substitute
$M=2k+2$ and get
$$
(k+3)k+1\geqslant 2k+3,
$$
which is true, either. Thus the claim of Proposition 2.3 is shown
for $i\leqslant M-2$.\vspace{0.1cm}

Finally, in the case $i=M-1$ the parameter $b$ takes the values
$0,1,\dots, M-3$ (it is precisely to treat the option $b=M-2$ that
we considered the case of a line separately in Proposition 2.2).
If $b=0$, we get $\varphi_1(0)=3M-3$ as above which is fine.
Substituting $b=M-3$, we get the value $\varphi_1(M-3)=M+3$ which
leads to the estimate
$$
\mathop{\rm codim}({\cal B}_{M-1}\subset{\cal A}) \geqslant M+1,
$$
as required. Now the proof of Proposition 2.3 is
complete.\vspace{0.1cm}

Q.E.D. for Theorem 0.2.


\begin{flushleft}
Department of Mathematical Sciences,\\
The University of Liverpool
\end{flushleft}

\noindent{\it pukh@liverpool.ac.uk}

\end{document}